\documentclass[12pt,oneside]{amsart}

\usepackage{amscd}
\usepackage{amssymb}
\usepackage[latin1]{inputenc}
\usepackage{verbatim}

\def\1{^{-1}}

\def\e{\varepsilon}
\def\f{\phi}

\def\G{\Gamma}

\def\cal{\mathcal}

\newtheorem{theorem}{Theorem}

\title{Wreath product of groups with infinite conjugacy classes}

\begin{document}
\maketitle
\begin{center}
{\sc Jean-Philippe PR\' EAUX}\footnote[1]{Centre de Recherche de
l'Armée de l'Air, Ecole de l'air, F-13661 Salon de Provence air}\
\footnote[2]{Laboratoire d'Analyse Topologie et Probabilités,
Universit\'e de Provence, 39 rue F.Joliot-Curie, F-13453 marseille
cedex 13\\
\indent {\it E-mail :} \ preaux@cmi.univ-mrs.fr}
\end{center}

\begin{abstract}
We characterise the group property of being with infinite
conjugacy
classes  for wreath products of groups.\\
\end{abstract}

\section{Introduction}

A group is said to be with {\sl infinite conjugacy classes} (or
{\sl icc}) if it is non trivial, and if all its conjugacy classes
except $\{ 1\}$ are infinite. This property is motivated by the
theory of Von Neumann algebra, since for any group $\G$, a
necessary and sufficient condition for its Von Neumann algebra to
be a type $II-1$ factor is that $\G$ be icc (cf. \cite{roiv}).
The property of being icc has been characterized in several
classes of groups : 3-manifolds and $PD(3)$ groups in
\cite{aogf3v}, and groups acting on Bass-Serre trees in
\cite{ydc}. We will focus here on groups defined as a wreath
product.

In the following, $D$, $Q$ are groups, $\Omega$ is a $Q$-set and
the group $G$ is  the wreath product $G=D\wr_\Omega Q$. That is,
let us denote by $D^{(\Omega)}$ the group of maps from $\Omega$ to
$D$ having a finite support and by $\lambda : Q\longrightarrow
Aut(D^{(\Omega)})$ the homomorphism  defined by $\forall x\in Q$,
$\forall\f\in D^{(\Omega)}$, $\lambda(q)(\f)(x)=\f(q^{-1}x)$ ; the
group $G$ is  defined to be the split extension
$G=D^{(\Omega)}\rtimes Q$ associated with $\lambda$, in the sense
that $\forall \f\in D^{(\Omega)}$, $\forall q\in Q$, $q\f
q^{-1}=\lambda(q)(\f)$.

\begin{theorem}
Let $G=D\wr_\Omega Q$, with $D\not=\{1\}$ ; a necessary and
sufficient condition for $G$  to be icc is that  on the one hand
condition (i) is satisfied :
\begin{itemize}
\item[\sl(i)]  1 is the only element of $FC(Q)$ which fixes $\Omega$
pointwise. \end{itemize}
 and on the other at least one of the
following conditions is satisfied :
\begin{itemize}
 \item[\sl(ii)] $D$ is icc,
 \item[\sl(iii)] all
$Q$-orbits in $\Omega$ are infinite.
\end{itemize}
In particular, if the action of $Q$ on $\Omega$ is free, then $G$
is icc if and only if either $D$ is icc or $Q$ is infinite.
\end{theorem}

I want to warmly acknowledge Pierre de la Harpe first for having
introduced me the problem, and second for all his useful advices
and
comments.\\

\section{Proof of the theorem}

Let us first give some notations : if $G$ is a group and $x,y$ are
element of $G$, then $x^y$ is the element of $G$ defined by
$x^y=y^{-1}xy$. If $H$ is a subgroup of $G$, then $x^H=\{x^y\ |\
y\in H\}$ ; in particular $x^G$ denote the conjugacy class of $x$
in $G$. The set of elements of $G$ having a finite conjugacy class
is a normal subgroup of $G$ that we denote $FC(G)$. \\

\noindent {\sl Proof.} In the following $\e$ will denote the
neutral element of $D^{(\Omega)}$, \emph{i.e.} the element of
$D^{(\Omega)}$ defined by $\forall x\in\Omega$, $\e(x)=1$. Given
$y\in \Omega$ and $d\in D$, the element $\zeta_d^y$ of
$D^{(\Omega)}$ is defined by $\zeta_d^y(x)=d$ if $x=y$ and
$\zeta_d^y(x)=1$ otherwise.


 We first suppose that $G$ is icc and prove the necessary part of the assumption.
 Necessarily condition $(i)$ is satisfied ; otherwise there would exist $q_0\not=1$ in
$FC(Q)$ fixing $\Omega$ pointwise, and $\{(\e,q_0^q)\in G\ |\ q\in
Q\}$ would be a finite subset of $G$ invariant under conjugacy,
contradicting that $G$ is icc. We now prove that if condition
$(iii)$ does not hold then condition $(ii)$ must hold.
Suppose  that  $\Omega$ has a finite $Q$-orbit $\cal O$. If $D$
would contain a finite conjugacy class $\xi\not=\{1\}$, then the
set $\Phi$ of maps from $\Omega$ to $\xi$ having their support in
$\cal O$ would be finite and non empty, and the subset
$\{(\f,1)\in G\ |\ \f\in\Phi\}$ of $G$ would be finite and
invariant under conjugacy, which is impossible. Hence if $G$ is
icc, either condition (ii) or condition (iii) is satisfied, which
proves the necessary part of the assumption.\smallskip\\
\indent We now prove the sufficient part of the assumption. So we
suppose in the following that condition $(i)$ is satisfied.

Suppose that condition $(iii)$ is satisfied, \emph{i.e.} all
$Q$-orbits in $\Omega$ are infinite. Let $g=(\f,q)\in G$ ; suppose
first that $\f\not=\e$. Its support is non empty and has an
infinite $Q$-orbit so that $\f\in D^{(\Omega)}$ has an infinity of
translated under the action of $\lambda(Q)$, and it follows that
$g^Q$ and hence also $g^G$ is infinite. Suppose now that
$g=(\e,q)$ is non trivial in $G$. If $q\not\in FC(Q)$ then $g^Q$
and hence also $g^G$ is infinite. If $q\in FC(Q)$ let $y\in\Omega$
be an element that $q$ does not fix (existence follows from
condition $(i)$),  and $d\not=1$ be an element of $D$.
Consider the element $g'=(\zeta_d^y,1)^{-1}g\
(\zeta_d^y,1)=(\phi,q)$ of $G$ ; $g'$ is a conjugate of $g$ and
$\phi=\zeta_{d^{-1}}^y\zeta_d^{qy}\not=\e$. Hence the above
argument applies to show that $g^G$ is infinite. It follows that
$G$ is icc.

Suppose that condition $(ii)$ is satisfied, \emph{i.e.} $D$ is
icc. Obviously each element $(\f,1)$ with $\f\not=\e$ has an
infinite conjugacy class. Let $g=(\f,q)$ with $q\not=1$. If
$q\not\in FC(Q)$ then $g^Q$ and $g^G$ are infinite. If $q\in
FC(Q)$, condition $(i)$ implies that $q$ does not fix some element
$y\in\Omega$.
%
Consider for any $d\in D$ the conjugate of $g$ that we denote by
$g_d=(\zeta_d^y,1)^{-1}(\f,q)\ (\zeta_d^y,1)$. If $y$ does not lie
in the support $Supp(\f)$ of $\f$, then
$g_d=(\f\zeta_{d^{-1}}^y\zeta_d^{qy},q)$. If $y\in Supp(\f)$, say
$\f=\f_0\zeta_c^y$ and $y\not\in Supp(\f_0)$, then
$g_d=(\f_0\zeta_{d^{-1}c}^y\zeta_d^{qy},q)$. In any case, since
$qy\not= y$ all $g_d$ for $d\in D$ are distinct. Since $D$ is icc
$D$ is infinite, and so $g$ has an infinite conjugacy class. Hence
$G$ is icc.
 \hfill$\square$

 \end{document}